\documentclass[a4paper]{article}

\pdfoutput=1

\usepackage[english]{babel}
\usepackage[utf8]{inputenc}
\usepackage[T1]{fontenc}

\usepackage{graphics}

\usepackage{amsmath}
\usepackage{amssymb}

\usepackage{tcolorbox}

\usepackage{authblk}

\title{Emulating complex networks with a single\\ delay differential equation}

\author[1,2]{Florian Stelzer\thanks{stelzer@math.tu-berlin.de}}

\author[1]{Serhiy Yanchuk}

\affil[1]{Institute of Mathematics, Technische Universit\"at Berlin}
\affil[2]{Department of Mathematics, Humboldt-Universit\"at zu Berlin}

\usepackage{cite}

\begin{document}

\maketitle

\begin{abstract}
	A single dynamical system with time-delayed feedback can emulate networks. This property of delay systems made them extremely useful tools for Machine Learning applications. 
	Here we describe several possible setups, which allow emulating multilayer (deep) feed-forward networks as well as recurrent networks of coupled discrete maps with arbitrary adjacency matrix by a single system with delayed feedback. While the network's size can be arbitrary, the generating delay system can have a low number of variables, including a scalar case. 
\end{abstract}

\section{Introduction}
\label{intro}

Systems with time-delays, or delay-differential equations (DDE), play an important role in modeling various natural phenomena and technological processes \cite{Stepan1989,Hale1993,Diekmann1995a,Erneux2009,Smith2010,Erneux2017,Yanchuk2017,Krisztin2008}. In optoelectronics, delays emerge due to finite optical or electric signal propagation time between the elements \cite{Vladimirov2005,Erzgraber2006,DHuys2008,Vicente2008,Fiedler2008,Wolfrum2010,Yanchuk2010a,Soriano2013,Oliver2015,Marconi2015,Puzyrev2016,Yanchuk2019}. Similarly, in neuroscience, propagation delays of the action potentials play a crucial role in information processing in the brain \cite{Foss2000,Wu2001,Izhikevich2006,Stepan2009,Deco2009a,Perlikowski2010a,Popovych2011,Kantner2013}. 

Machine Learning is another rapidly developing application area of delay systems  \cite{Paugam-Moisy2008,Appeltant2011,Martinenghi2012,Appeltant2012a,Larger2012,Brunner2013,SCH13l,Toutounji2014,Grigoryeva2015,Penkovsky2017,Larger2017,Harkhoe2019,Stelzer2019,Hart2019,Koster2020,Koester2020,Argyris2020,Goldmann2020,Sugano2020}. 
It is shown recently that DDEs can successfully realize a reservoir computing setup, theoretically \cite{Hart2017,Keuninckx2017,Hart2019,Stelzer2019,Koster2020,Stelzer2020,Koester2020,Goldmann2020}, and implemented in optoelectronic hardware \cite{Appeltant2011,Appeltant2012a,Larger2017}.
In time-delay reservoir computing, a single DDE with either one or a few variables is used for building a ring network of coupled maps with fixed internal weights and fixed input weights. 
In a certain sense, the network structure emerges by properly unfolding the temporal behavior of the DDE. 
In this paper, we explain how such an unfolding appears, not only for the ring network as in  reservoir computing but also for arbitrary networks of coupled maps.
In \cite{Hermans2015}, a training method is proposed to modify the input weights while the internal weights are still fixed.

Among the most related previous publications, Hart and collaborators unfold networks with arbitrary topology from delay systems \cite{Hart2017,Hart2019}.
Our work extends their results in several directions, including varying coupling weights and applying it to a broader class of delay systems.
The networks constructed by our method allow for a modulation of weights. Hence,  they can be employed in Machine Learning applications with weight training. 
In our recent paper \cite{Stelzer2020}, we show that a single DDE can emulate a deep neural network and perform various computational tasks successfully. More specifically, the work \cite{Stelzer2020} derives a multilayer neural network from a delay system with modulated feedback terms. This neural network is trained by gradient descent using back-propagation and applied to machine learning tasks.

As follows from the above-mentioned machine learning applications, delay models can be effectively used for unfolding complex network structures in time. 
Our goal here is a general description of such networks. 
While focusing on the network construction, we do not discuss details of specific machine learning applications such as e.g., weights training by gradient descent, or specific tasks. 

The structure of the paper is as follows. 
In Sec.~\ref{sec:general-case} we derive a feed-forward network from a DDE with modulated feedback terms. Section~\ref{sec:recurrent-network} describes a recurrent neural network. 
In Sec.~\ref{sec:semilinear-case}, we review a special but practically important case of delay systems with a linear instantaneous part and nonlinear delayed feedback containing an affine combination of the delayed variables; originally, these results have been derived in \cite{Stelzer2020}.

\section{From delay systems to multilayer feed-forward networks}
\label{sec:general-case}

\subsection{Delay systems with modulated feedback terms}
\label{subsec:delay-system}

Multiple delays are required for the construction of a network with arbitrary topology by a delay system \cite{Hart2017,Hart2019,Stelzer2020}. In such a network, the connection weights are emulated by a modulation of the delayed feedback signals \cite{Stelzer2020}.
Therefore, we consider a DDE of the following form 
\begin{align}\label{eq:general-system}
	\dot{x}(t) = f(x(t), z(t), \mathcal{M}_1(t)x(t-\tau_1),\ldots , \mathcal{M}_D(t)x(t-\tau_D)),
\end{align}
with $D$ delays $\tau_1,\dots,\tau_D$, a nonlinear function $f$, a time-dependent driving signal $z(t)$, and modulation functions $\mathcal{M}_1(t), \ldots , \mathcal{M}_D(t)$. 

\begin{figure}
    \centering
    \includegraphics{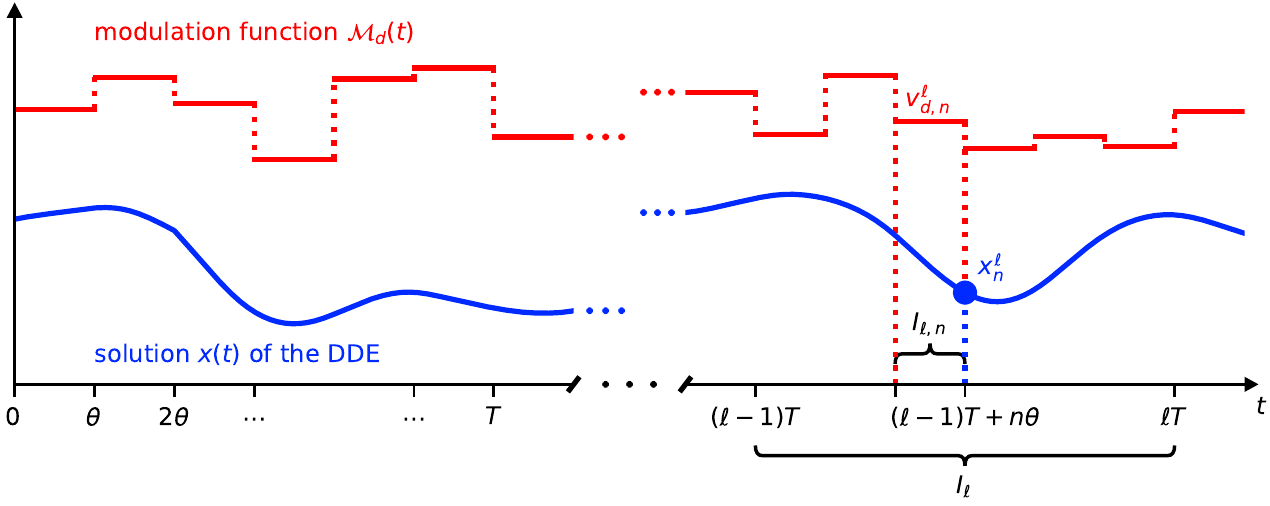}
    \caption{
    Illustration of the clock cycle intervals $I_\ell$ and  sub-intervals $I_{\ell,n}$. The node $x^\ell_n$ (blue dot) is defined by the value of the solution $x(t)$ of system~\eqref{eq:general-system} (blue line) at the time point $t=(\ell -1)T+n\theta$. The modulation function $\mathcal{M}_d(t)$ is a step function with constant values $v^\ell_{d,n}$ on the intervals $I_{\ell ,n}$.
    }
    \label{fig:first}
\end{figure}

System \eqref{eq:general-system} is a non-autonomous DDE, and the properties of the functions $\mathcal{M}_d(t)$ and $z(t)$ play an important role for unfolding a network from \eqref{eq:general-system}. 
To define these properties, a time quantity $T>0$ is introduced, called the \textit{clock cycle}. 
Further, we choose a number $N$ of grid points per $T$-interval and define $\theta := T/N$. 
We define the clock cycle intervals 
$$I_\ell := ((\ell-1)T, \ell T],\ \ell = 1, \ldots ,L,$$ 
which we split into smaller sub-intervals 
$$I_{\ell, n} := ((\ell -1)T + (n-1)\theta, (\ell - 1)T + n\theta], \ n=1,\ldots ,N,$$
see Fig.~\ref{fig:first}.
We assume the following properties for the delays and modulation functions:\\
\textbf{Property (I):} The delays satisfy $\tau_d = n_d \theta,\ d=1,\ldots ,D$ with natural numbers $0 < n_1 < \cdots < n_D < 2N$. Consequently, it holds $0 < \tau_1 < \cdots < \tau_D < 2T$. \\
\textbf{Property (II):} The functions $\mathcal{M}_d(t)$ are step-functions, which are constant on the intervals $I_{\ell,n}$. We denote these constants as $v^\ell_{d,n}$, i.e.
$$
\mathcal{M}_d(t) = v^\ell_{d,n} \quad \text{for} \quad t\in I_{\ell,n}.
$$

In the following sections, we show that one can consider the intervals $I_\ell$ as layers with $N$ nodes of a network arising from the delay system~\eqref{eq:general-system} if the modulation functions $\mathcal{M}_d(t)$ fulfill certain additional requirements. The $n$-th node of the $\ell$-th layer is defined as
\begin{align}
\label{eq:nodes}
	x^\ell_n := x((\ell - 1)T + n\theta), \quad n = 1,\ldots ,N, \ \ell =1,\ldots ,L,
\end{align}
which corresponds to the solution of the DDE \eqref{eq:general-system} at time point $(\ell - 1)T + n\theta$. 
The solution at later time points $x^{\ell'}_{n'}$ with either $\ell'>\ell$ or $n'>n$ for $\ell'=\ell$ depends, in general, on $x^\ell_n$, thus, providing the interdependence between the nodes.
Such dependence can be found explicitly in some situations. 
The simplest way is to use a discretization for small $\theta$, and we consider such a case in the following Sec.~\ref{subsec:discretization}. Another case, when $\theta$ is large, can be found in \cite{Stelzer2020}.

Let us remark about the initial state for DDE~\eqref{eq:general-system}. According to the general theory~\cite{Hale1993}, in order to solve an initial value problem, an initial history function $x_0(s)$ must be provided on the interval  $s \in [-\tau_D, 0]$, where $\tau_D$ is the maximal delay. 
In terms of the nodes, one needs to specify $x_n^\ell$ for $n_D$ ``history'' nodes.
However, the modulation functions $\mathcal{M}_d(t)$ can weaken this requirement. For example, if 
$\mathcal{M}_d(t) = 0$ for $t \leq \tau_d$, then it is sufficient to know the initial state $x(0) = x_0^1 =x_0$ at a single point, and we do not require a history function at all. In fact, the latter special case has been employed in \cite{Stelzer2020} for various machine learning tasks.

\subsection{Disclosing network connections via discretization of the DDE}
\label{subsec:discretization}

Here we consider how a network of coupled maps can be derived from  DDE \eqref{eq:general-system}. Since the network nodes are already introduced in Sec.~\ref{subsec:delay-system} as $x_n^\ell$ by Eq.~\eqref{eq:nodes}, it remains to describe the connections between the nodes. Such links are functional connections between the nodes $x_n^\ell$. Hence, our task is to find functional relations (maps) between the nodes. 

For simplicity, we restrict ourselves to the Euler discretization scheme since 
the obtained network topology is independent of the chosen discretization.
Similar network constructions by discretization from ordinary differential equations have been employed in \cite{Haber2017,Lu2018,Chen2018}. 

We apply a combination of the forward and backward Euler method: the instantaneous system states of~\eqref{eq:general-system} are approximated by the left endpoints of the small-step intervals of length $\theta$ (forward scheme). The driving signal $z(t)$ and the delayed system states are approximated by the right endpoints of the step intervals (backward scheme). Such an approach leads to simpler expressions. We obtain
\begin{align}\label{eq:euler}
	x^\ell_n = x^\ell_{n-1} + \theta f(x^\ell_{n-1},z(t^\ell_n), \mathcal{M}_1(t^\ell_n)x(t^\ell_n-\tau_1),\ldots ,\mathcal{M}_D(t^\ell_n)x(t^\ell_n-\tau_D))
\end{align}
for $n = 2,\ldots,N$, where $t^\ell_n := (\ell - 1)T + n\theta$, and
\begin{align}\label{eq:euler1}
x^\ell_1 = x^{\ell-1}_N + \theta f(x^{\ell-1}_N, z(t^\ell_1), \mathcal{M}_1(t^\ell_1)x(t^\ell_1-\tau_1),\ldots ,\mathcal{M}_D(t^\ell_1)x(t^\ell_1-\tau_D))
\end{align}
for the first node in the $I_\ell$-interval.

According to Property (I), the delays satisfy $0<\tau_d<2T$. Therefore, the delay-induced feedback connections with target in the interval $I_\ell$ can originate from one of the following intervals: $I_\ell$, $I_{\ell-1}$, or $I_{\ell-2}$. In other words: the time points $t^\ell_n - \tau_d$ can belong to one of these intervals $I_{\ell}$, $I_{\ell -1}$, $I_{\ell -2}$. 
Formally, it can be written as
\begin{equation}
\label{eq:inIell}
t_{n}^{\ell}-\tau_{d}=t_{n}^{\ell}-n_{d}\theta=\begin{cases}
t_{n-n_{d}}^{\ell}\in I_{\ell}, & \text{if}\quad n_{d}<n,\\
t_{N+n-n_{d}}^{\ell-1}\in I_{\ell-1}, & \text{if}\quad n\leq n_{d}<N+n,\\
t_{2N+n-n_{d}}^{\ell-2}\in I_{\ell-2}, & \text{if}\quad N+n\leq n_{d}.
\end{cases}
\end{equation}

We limit the class of networks to multilayer systems with connections between the neighboring layers. Such networks, see Fig.~\ref{fig:network-from-delay-system}b, are frequently employed in machine learning tasks, e.g. as deep neural networks \cite{Bishop2006,Goodfellow2016,Lecun2015,Schmidhuber2015,Stelzer2020}.
Using \eqref{eq:inIell}, we can formulate a condition for the modulation functions $\mathcal{M}_d(t)$ to ensure that the delay terms $x(t-\tau_d)$ induce only connections between subsequent layers. 
For this, we set the modulation functions' values to zero if the originating time point $t^\ell_n - \tau_d$ of the corresponding delay connection does not belong to the interval $I_{\ell-1}$. 
This leads to the following assumption on the modulation functions:\\
\textbf{Property (III):} The modulation functions $\mathcal{M}_d(t)$ vanish at the following intervals:
\begin{align}\label{eq:M-condition}
	\mathcal{M}_d(t) = v_{d,n}^\ell = 0 \quad \text{for} \quad t\in I_{\ell,n} \quad \text{if} \quad (n_d < n) \quad \text{or} \quad (N+n \leq n_d).
\end{align}
In the following, we assume that condition (III) is satisfied.

Expressions \eqref{eq:euler}--\eqref{eq:euler1} contain the interdependencies between $x_n^\ell$, i.e., the connections between the nodes of the network. We explain these dependencies and present them in a more explicit form in the following. Our goal is to obtain the multilayer network shown in Fig.~\ref{fig:network-from-delay-system}b. 

\subsection{Effect of time-delays on the network topology}
\label{subsec:network-topology}

\begin{figure}
	\centering
	\includegraphics{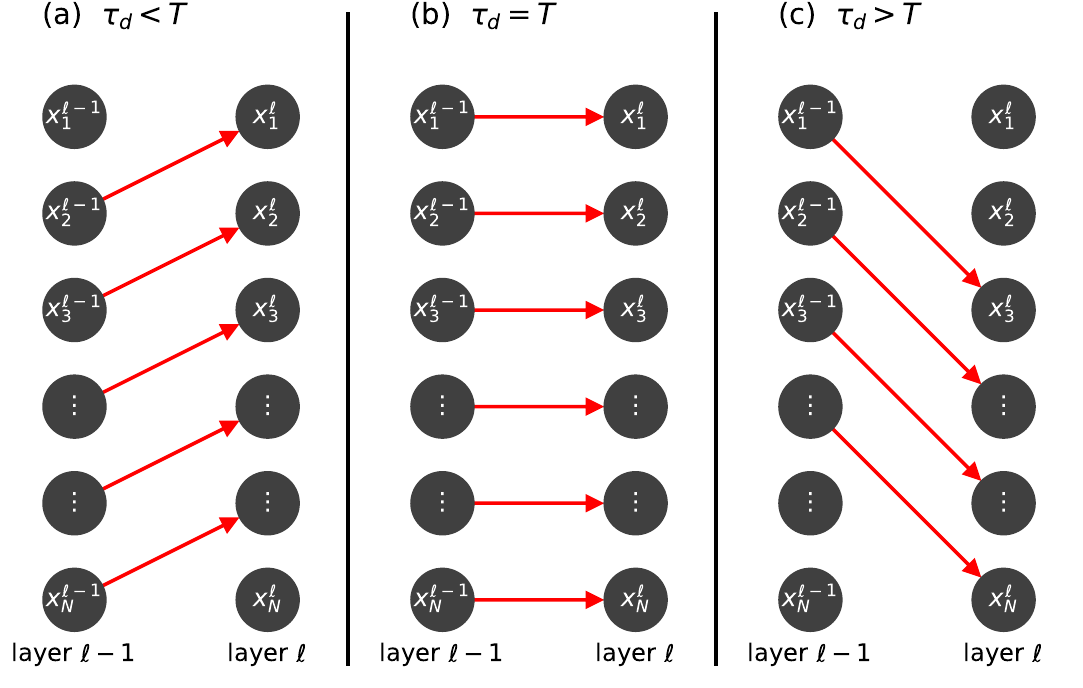}
	\caption{Network connections induced by one time-delay $\tau_d$. Panel (a): connections induced by $\tau_d<T$. Panel (b):  $\tau_d = T$. Panel (c): $\tau_d > T$. Multiple delays $\tau_1, \ldots ,\tau_D$ result in a superposition of parallel patterns as shown in Fig.~\ref{fig:network-from-delay-system}b.
	}
	\label{fig:net-topo}       
\end{figure}

Taking into account property (III), the node $x^\ell_n$ of layer $I_{\ell}$ receives a connection from a node $x^{\ell -1}_{n-n'_d}$ of layer $I_{\ell -1}$, where $n'_d:=n_d-N$. Two neighboring layers are illustrated in Fig.~\ref{fig:net-topo}, where the nodes in each layer are 
ordered vertically from top to bottom. 
Depending on the size of the delay, we can distinguish three cases.
\begin{itemize}
	\item[(a)] For $\tau_d<T$, there are $n_d$ ``upward'' connections as shown in panel Fig.~\ref{fig:net-topo}a. 
	\item[(b)] For $\tau_d=T$, there are $n_d=N$  ``horizontal'' delay-induced connections, i.e. connections from nodes of layer $\ell -1$ to nodes of layer $\ell$ with the same index, see Fig.~\ref{fig:net-topo}b.
	\item[(c)] For larger delays $\tau_d>T$, there are $2N-n_d$ ``downward'' delay-induced connections, as shown in Fig.~\ref{fig:net-topo}c.
\end{itemize}
In all cases, the connections induces by one delay $\tau_d$ are parallel. Since the delay system possesses multiple delays $0 < \tau_1 < \ldots < \tau_D < 2T$, the parallel connection patterns overlap, as illustrated in Fig.~\ref{fig:network-from-delay-system}b, leading to a more complex topology. In particular, a fully connected pattern appears for $D = 2N-1$ and $\tau_d = \theta d$.

\subsection{Modulation of connection weights}
\label{subsec:weight-modulation}

With the modulation functions satisfying property (III), the Euler scheme~\eqref{eq:euler}--\eqref{eq:euler1} simplifies to the following map
\begin{align}\label{eq:general-map1}
x^\ell_1 &= x^{\ell-1}_N + \theta f(x^{\ell-1}_N,z(t^\ell_1), v^\ell_{1,1} x^{\ell -1}_{1-n'_1},\ldots ,v^\ell_{D,1} x^{\ell -1}_{1-n'_D}),\\
x^\ell_n &= x^\ell_{n-1} + \theta f(x^\ell_{n-1}, z(t^\ell_n), v^\ell_{1,n} x^{\ell -1}_{n-n'_1},\ldots ,v^\ell_{D,n} x^{\ell -1}_{n-n'_D}), \quad n=2,\ldots ,N,\label{eq:general-map}
\end{align}
where Eq.~\eqref{eq:M-condition} implies $v^{\ell}_{d,n} = 0$ if $n-n'_d < 1$ or $n-n'_d > N$. In other words, the dependencies at the right-hand side of \eqref{eq:general-map1}--\eqref{eq:general-map} contain only the nodes from the $\ell-1$-th layer.  Moreover, the numbers $v^\ell_{d,n}$ determine the strengths of the connections from $x^{\ell -1}_{n-n'_d}$ to  $x^\ell_n$ and can be considered as network weights. By reindexing, we can define weights $w^\ell_{nj}$ connecting node $j$ of layer $\ell -1$ to node $n$ of layer $\ell$. These weights are given by the equation
\begin{align}\label{eq:weight-matrix1}
w^\ell_{nj} := \sum_{d =1}^D \delta_{n-n'_d,j} v^\ell_{d,n} 
=
\begin{cases}
0 & \text{if } \forall d \colon j \neq n - n'_d,\\
v^\ell_{d,n} & \text{if } \exists d \colon j = n - n'_d,
\end{cases}
\end{align}
and define the entries of the weight matrix $W^\ell = (w^\ell_{nj}) \in \mathbb{R}^{N\times (N+1)}$, except for the last column, which is defined below and contains bias weights. The symbol $\delta_{nj}$ is the Kronecker delta, i.e. $\delta_{nj} = 1$ if $n=j$, and $\delta_{nj} = 0$ if $n\neq j$.

\begin{figure}[t]
	\centering
	\includegraphics{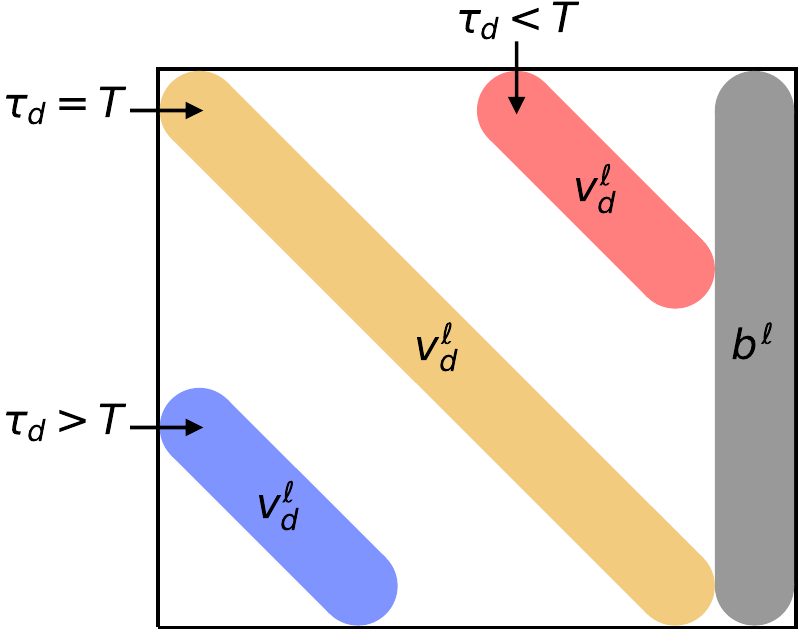}
	\caption{
	Coupling matrix $W^\ell$ between the hidden layers $\ell - 1$ and $\ell$, see Eq.~\eqref{eq:weight-matrix1}--\eqref{eq:weight-matrix2}.  The nonzero weights are arranged along the diagonals, and equal $v^\ell_{d,n}$. The position of the diagonals is determined by the corresponding delay $\tau_d$. 
	If $\tau_d = T = N\theta$, then the main diagonal contains the entries $v^\ell_{d,1},\ldots ,v^\ell_{d,N}$ (shown in yellow). If $\tau_d = n_d\theta < T$, then the corresponding diagonal lies above the main diagonal and contains the values $v^\ell_{d , 1}, \ldots , v^\ell_{d , n_d}$ (red). If $\tau_d = n_d\theta > T$, then the corresponding diagonal lies below the main diagonal and contains the values $v^\ell_{d , n_d-N+1}, \ldots , v^\ell_{d, N}$ (blue). The last column of the matrix contains the bias weights (gray).
	}
	\label{fig:weight-matrix}
\end{figure}

The time-dependent driving function $z(t)$ can be utilized to realize a bias weight $b^\ell_n$ for each node $x^\ell_n$. For details, we refer to  Sec.~\ref{subsec:mulilayer-network}. We define the last column of the weight matrix $W^\ell$ by
\begin{align}\label{eq:weight-matrix2}
w^\ell_{n,N+1} := b^\ell_n.
\end{align}
The weight matrix is illustrated in Fig.~\ref{fig:weight-matrix}.
This matrix $W^\ell$ is in general sparse, where the degree of sparsity depends on the number $D$ of delays. If $D=2N-1$ and $\tau_d = d\theta, \ d=1,\ldots ,D$, we obtain a dense connection matrix. Moreover, the positions of the nonzero entries and zero entries are the same for all matrices $W^2, \ldots , W^L$, but the values of the nonzero entries are in general different. 

\subsection{Interpretation as multilayer neural network}
\label{subsec:mulilayer-network}

The map~\eqref{eq:general-map1}--\eqref{eq:general-map} can be interpreted as the hidden layer part of a multilayer neural network provided we define suitable input and output layers.

\begin{figure}[t]
	\centering
	\includegraphics{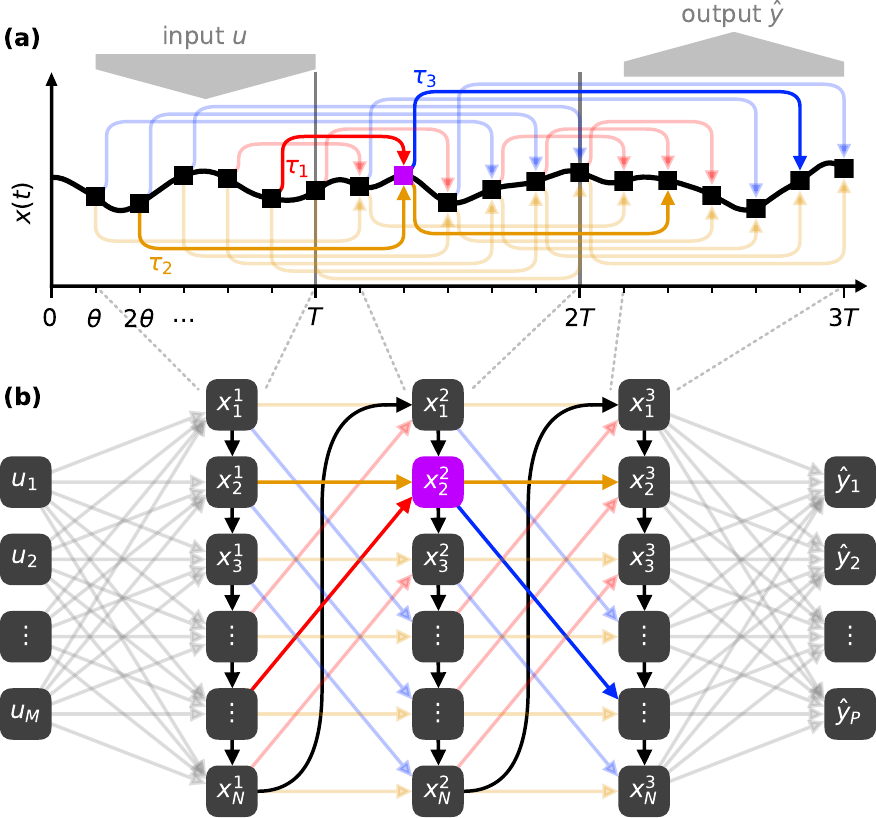}
	\caption{Implementing a multilayer neural network by delay system~\eqref{eq:general-system}. 		Panel (a): The system state is considered at discrete time points $x^\ell_{n} := x((\ell - 1)T +n\theta)$. The intervals $I_\ell$ correspond to layers. Due to delayed feedback, non-local connections emerge (color lines). Panel (b) shows the resulting neural network.
	}
	\label{fig:network-from-delay-system}
\end{figure}

The input layer determines how a given input vector $u\in \mathbb{R}^{M+1}$ is transformed to the state of the first hidden layer $x(t),\ t\in I_1$.
The input $u\in \mathbb{R}^{M+1}$
contains $M$ input values $u_1,\ldots ,u_M$ and an additional entry $u_{M+1}=1$.  
In order to ensure that $x(t),\ t\in I_1$ depends on $u$ and the initial state $x(0)=x_0$ exclusively, and does not depend on a history function $x(s),\ s < 0$, we set all modulation functions to zero on the first hidden layer interval. This leads to the following\\
\textbf{Property (IV):}
The modulation functions satisfy 
\begin{align}
	\mathcal{M}_d(t) = 0, \quad t\in I_1, \ d=1,\ldots ,D.
\end{align}
The dependence on the input vector $u\in\mathbb{R}^{M+1}$ can be realized by the driving signal $z(t)$. \\
\textbf{Property (V):} The driving signal $z(t)$ on the interval $I_1$ is the step function given by
\begin{align}\label{eq:J-1}
z(t) = & J(t) \quad  \text{for } \quad t\in I_1,     \\
& J(t) = J_n = \left[ f^\mathrm{in}(W^{\mathrm{in}} u) \right]_n \quad \text{for} \quad t\in I_{1,n},\label{eq:J-2}
\end{align}
where $f^\mathrm{in}(W^{\mathrm{in}} u)\in\mathbb{R}^N$ is the preprocessed input,  $W^{\mathrm{in}}\in\mathbb{R}^{N\times (M+1)}$ is an input weight matrix, and $f^\mathrm{in}$ is an element-wise input preprocessing function. For example,  $f^\mathrm{in}(a)=\tanh (a)$ was used in~\cite{Stelzer2020}. 

As a result, the following holds for the first hidden layer
\begin{align}\label{eq:general-system-first-layer}
\dot{x}(t) = f(x(t), J(t), 0,\ldots ,0), \quad t\in I_1,
\end{align}
which is just a system of ordinary differential equations, which requires an initial condition at a single point $x(0)=x_0$ for solving it in positive time. 
This yields the coupled map representation
\begin{align}\label{eq:general-map1-first-layer}
x^\ell_1 &= x_0 + \theta f(x_0, J_{1},0,\ldots ,0),\\
x^1_n &= x^1_{n-1} + \theta f(x^1_{n-1}, J_{n},0,\ldots ,0), \quad n=2,\ldots
,N.\label{eq:general-map-first-layer}
\end{align}

For the hidden layers $I_2,I_3,\dots$, the driving function $z(t)$ can be used to introduce a bias as follows. \\
\textbf{Property (VI):} The  driving  signal $z(t)$  on  the  intervals $I_\ell$, $\ell\ge 2$,  is  the  step  function given by
\begin{align}
z(t) = & b(t) \quad  \text{for } \quad t>T,     \\
& b(t) = b_n^\ell \quad \text{for} \quad t\in I_{\ell,n},\quad  \ell \ge 2.
\end{align}
Assuming the properties (I)--(VI), Eqs.~\eqref{eq:general-map1}--\eqref{eq:general-map} imply
\begin{align}\label{eq:general-map1-hidden-layer}
x^\ell_1 &= x^{\ell-1}_N + \theta f(x^{\ell-1}_N,b^\ell_1, v^\ell_{1,1} x^{\ell -1}_{1-n'_1},\ldots ,v^\ell_{D,1} x^{\ell -1}_{1-n'_D}),\\
x^\ell_n &= x^\ell_{n-1} + \theta f(x^\ell_{n-1}, b^\ell_n, v^\ell_{1,n} x^{\ell -1}_{n-n'_1},\ldots ,v^\ell_{D,n} x^{\ell -1}_{n-n'_D}), \quad n=2,\ldots ,N.\label{eq:general-map-hidden-layer}
\end{align}

Let us finally define the output layer, which transforms the node states $x^L_1,\ldots , x^L_n$ of the last hidden layer to an output vector $\hat{y} \in \mathbb{R}^P$. For this, we define a vector $x^L := (x^L_1, \ldots , x^L_N, 1)^\mathrm{T} \in \mathbb{R}^{N+1}$, an output weight matrix $W^\mathrm{out}\in\mathbb{R}^{P\times (N+1)}$, and an output activation function $f^\mathrm{out}\colon \mathbb{R}^P \to \mathbb{R}^P$. The output vector is then defined as
\begin{align}\label{eq:general-output}
	\hat{y} = f^\mathrm{out}(W^\mathrm{out}x^L).
\end{align}

Figure~\ref{fig:network-from-delay-system} illustrates the whole construction process of the coupled maps network; it is given by the equations~\eqref{eq:general-map1-first-layer}--\eqref{eq:general-output}.

We summarize the main result of section 2.

\begin{tcolorbox}
Under assumptions (I)--(VI) and for small $\theta$, DDE \eqref{eq:general-system} describes the multilayer network of coupled maps shown in Fig.~\ref{fig:network-from-delay-system}, with the specific dependencies given by Eqs.~\eqref{eq:general-map1-first-layer}, \eqref{eq:general-map-first-layer}, \eqref{eq:general-map1-hidden-layer}, \eqref{eq:general-map-hidden-layer},  and \eqref{eq:general-output}.
\end{tcolorbox}

\section{Constructing a recurrent neural network from a delay system}
\label{sec:recurrent-network}

System~\eqref{eq:general-system} can also be considered as recurrent neural network. To show this, we consider the system on the time interval $[0,KT]$, for some $K\in\mathbb{N}$, which is divided into intervals $I_k :=((k-1)T,kT], \ k = 1, \ldots ,K$. 
We use $k$ instead of $\ell$ as index for the intervals to make clear that the intervals do not represent layers. The state $x(t)$ on an interval $I_k$ is interpreted as the state of the recurrent network at time $k$. More specifically, 
\begin{align}
x^k_n := x((k - 1)T + n\theta), \quad n = 1,\ldots ,N, \ k = 1, \ldots ,K
\end{align}
is the state of node $n$ at the discrete time $k$. 
The driving function $z(t)$ can be utilized as an input signal for each $k$-time-step.\\
\textbf{Property (VII):}
$z(t)$ is the $\theta$-step function with 
\begin{align}
z(t)= z_n^k & \quad \text{for} \quad t\in I_{k,n}, \\
& (z_1^k,\dots,z_N^k)^T = f^\mathrm{in}(W^\mathrm{in} u(k)),
\end{align}
where $u(k),\ k = 1, \ldots ,K$ are $(M+1)$-dimensional input vectors, $W^\mathrm{in} \in \mathbb{R}^{N\times (M+1)}$ is an input weight matrix, $f^\mathrm{in}$ is an element-wise input preprocessing function. Each input vector $u(k)$ contains $M$ input values $u_1(k),\ldots ,u_M(k)$ and a fixed entry $u_{M+1}(k):=1$ which is needed to include bias weights in the last column of $W^\mathrm{in}$.

The main difference of the Property (VII) from (VI) is that it allows for the information input through $z(t)$ in all intervals $I_k$. Another important difference is related to the modulation functions, which must be $T$-periodic in order to implement a recurrent network. This leads to the following assumption. \\
\textbf{Property (VIII):}
The modulation functions $\mathcal{M}_d(t)$ are $T$-periodic $\theta$-step functions with 
\begin{align}
\mathcal{M}_d(t)= v_{d,n}\quad \text{for} \quad t\in I_{k,n}. 
\end{align}
Note that the value $v_{d,n}$ is independent on $k$ due to periodicity of $\mathcal{M}_d(t)$. When assuming the Properties (I), (III), (IV), (VII), and (VIII), the map equations~\eqref{eq:general-map1}--\eqref{eq:general-map} become
\begin{align}\label{eq:general-map1-recurrent}
x^k_1 &= x^{k-1}_N + \theta f(x^{k-1}_N, z^k_1, v_{1,1} x^{k -1}_{1+n'_1},\ldots ,v_{D,1} x^{k -1}_{1+n'_D}),\\
x^k_n &= x^k_{n-1} + \theta f(x^k_{n-1}, z^k_n, v_{1,n} x^{k -1}_{n+n'_1},\ldots ,v_{D,n} x^{k -1}_{n+n'_D}), \quad n=2,\ldots ,N,\label{eq:general-map-recurrent}
\end{align}
and can be interpreted as a recurrent neural network with the input matrix $W^\mathrm{in}$ and the internal weight matrix $W = (w_{nj}) \in \mathbb{R}^{N\times N}$ defined by
\begin{align}\label{eq:recurrent-weight-matrix}
w_{nj} := \sum_{d =1}^D \delta_{n-n'_d,j} v_{d,n}
=
\begin{cases}
0 & \text{if } \forall d \colon j \neq n - n'_d,\\
v_{d,n} & \text{if } \exists d \colon j = n - n'_d.
\end{cases}
\end{align}

\begin{figure}
	\centering
	\includegraphics{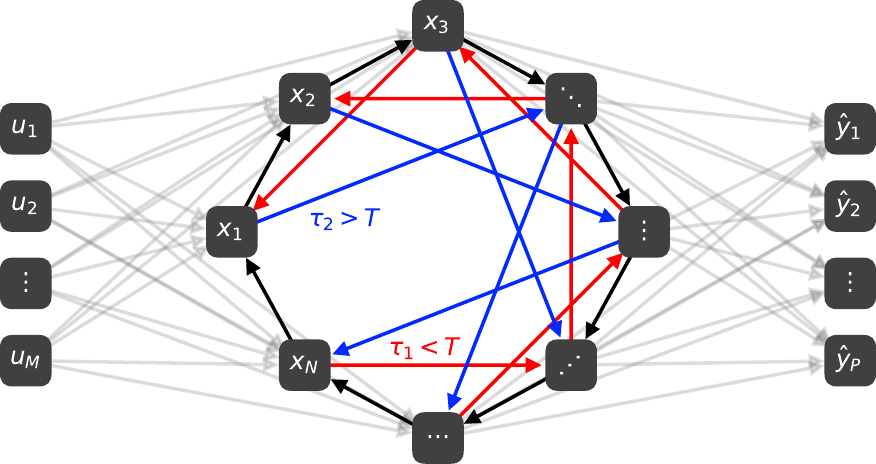}
	\caption{Recurrent network obtained from DDE \eqref{eq:general-system} with two delays. The delays $\tau_1<T$ and $\tau_2>T$ induce connections with opposite direction (color arrows). Moreover, the nodes of the recurrent layer a linearly locally coupled (black arrows). All nodes of the recurrent layers are connected to the input and output layer.
	}
	\label{fig:recurrent-network}
\end{figure}

When we choose the number of delays to be $D=2N-1$, we can realize any given connection matrix $W\in \mathbb{R}^{N\times N}$. For that we need to choose the delays $\tau_d = d\theta , \ d=1,\ldots ,2N-1$. Consequently there are $D=2N-1$ modulation functions $\mathcal{M}_d(t)$ which are step functions with values $v_{d,n}$. In this case Eq.~\eqref{eq:recurrent-weight-matrix} provides for all entries $w_{nj}$ of $W$ exactly one corresponding $v_{d,n}$. Therefore, the arbitrary matrix $W$ can be realized by choosing appropriate step heights for the modulation functions. In the setting of Sec.~\ref{sec:recurrent-network} the resulting network is an arbitrary recurrent network.

Summarizing, the main message of Sec. 3 is as follows.

\begin{tcolorbox}
Under assumptions (I), (III), (IV), (VII), and (VIII), and for small $\theta$, DDE \eqref{eq:general-system} describes the recurrent network shown in Fig.~\ref{fig:recurrent-network}, with the specific dependencies given by Eqs.~\eqref{eq:general-map1-recurrent}--\eqref{eq:general-map-recurrent} and an internal weight matrix $W$ given by \eqref{eq:recurrent-weight-matrix}.
\end{tcolorbox}

\section{Networks from delay systems with linear instantaneous part and nonlinear delayed feedback}
\label{sec:semilinear-case}
\label{subsec:semilinear-delay-system}

Particularly suitable for the construction of neural networks are delay systems with a stable linear instantaneous part and a feedback given by a nonlinear function of an affine combination of the delay terms and a driving signal. Such DDEs are described by the equation
\begin{align}\label{eq:semilinear-system}
\dot{x}(t) &= -\alpha x(t) + f(a(t)),
\end{align}
where $\alpha > 0$ is a constant time scale, $f$ is a nonlinear function, and
\begin{align}\label{eq:activation-signal}
a(t) &= z(t) + \sum_{d=1}^D \mathcal{M}_d(t)x(t-\tau_d).
\end{align}
Ref.~\cite{Krisztin2008} studied this type of equation for the case $D=1$, i.e. for one delay.

An example of \eqref{eq:semilinear-system} is the Ikeda system~\cite{Ikeda1979} where $D=1$, i.e. $a(t)$ consists of only one scaled feedback term $x(t-\tau)$, signal $z(t)$, and the nonlinear function $f(a)=\sin(a)$. This type of dynamics can be applied to reservoir computing using optoelectronic hardware~\cite{Larger2012}.
Another delay dynamical system of type~\eqref{eq:semilinear-system}, which can be used for reservoir computing, is the Mackey-Glass system~\cite{Appeltant2011}, where $D=1$ and the nonlinearity is given by $f(a)=\eta a /(1+|a|^p)$ with constants $\eta , p > 0$.
In the work~\cite{Stelzer2020}, system~\eqref{eq:semilinear-system} is used to implement a deep neural network. 

Even though the results of the previous sections are applicable to \eqref{eq:semilinear-system}--\eqref{eq:activation-signal}, the special form  of these equations allows for an alternative, more precise approximation of the network dynamics.

\subsection{Interpretation as multilayer neural network}
\label{subsec:ddnn-network}

It is shown in \cite{Stelzer2020} that one can derive a particularly simple map representation for system~\eqref{eq:semilinear-system} with activation signal~\eqref{eq:activation-signal}. We do not repeat here the derivation, and only present the resulting expressions. By applying a semi-analytic Euler discretization and the variation of constants formula, the following  equations connecting the nodes in the network are obtained:
\begin{align}\label{eq:first-hidden-layer-first-node}
x^1_1 &= e^{-\alpha \theta} x_0 + \alpha^{-1}(1-e^{-\alpha \theta}) f( a^1_1), \\
x^1_n &= e^{-\alpha \theta} x^1_{n-1} + \alpha^{-1}(1-e^{-\alpha \theta}) f( a^1_n), \quad n = 2, \ldots , N,
\end{align}
for the first hidden layer. 
The hidden layers $\ell = 2, \ldots ,L$ are given by
\begin{align}\label{eq:hidden-layer-first-node}
x^\ell_1 &= e^{-\alpha \theta} x^{\ell -1}_N + \alpha^{-1}(1-e^{-\alpha \theta}) f( a^\ell_1), \\
x^\ell_n &= e^{-\alpha \theta} x^\ell_{n-1} + \alpha^{-1}(1-e^{-\alpha \theta}) f( a^\ell_n), \quad n = 2, \ldots , N.\label{eq:hidden-layer}
\end{align}	
The output layer is defined by
\begin{align}\label{eq:output-layer}
\hat{y}_p := f^\mathrm{out}_p (a^\mathrm{out}), \quad p=1,\ldots, P,
\end{align}
where $f^\mathrm{out}$ is an output activation function.
Moreover,
\begin{align}
a^\mathrm{in}_n &:= \sum_{m=1}^{M+1} w^\mathrm{in}_{nm} u_m, & & n = 1, \ldots , N, \label{eq:activation-in}\\
a^1_n &:= g(a^\mathrm{in}_n), & & n = 1, \ldots , N, \label{eq:activation-first}\\
a^\ell_n &:= \sum_{j=1}^{N+1} w^\ell_{nj} x^{\ell -1}_j, & & n=1, \ldots , N, \ \ell = 2, \ldots , L,\label{eq:activation}\\
a^\mathrm{out}_p &:= \sum_{n=1}^{N+1} w^\mathrm{out}_{pn} x^L_n,  & & p = 1, \ldots , P,\label{eq:activation-output}
\end{align}
where $u_{M+1}:=1$ and $x^\ell_{N+1} := 1$, for $\ell=1,\ldots ,L$.

One can also formulate the relation between the hidden layers in a matrix form. For this, we define
\begin{align}
	A := \begin{pmatrix}
	0 & \cdots & \cdots & \cdots & 0 \\
	e^{-\alpha\theta} & \ddots &  & & \vdots \\
	0 & \ddots & \ddots &  & \vdots \\
	\vdots & \ddots & \ddots & \ddots &\vdots  \\
	0 & \cdots & 0 & e^{-\alpha\theta} & 0
	\end{pmatrix}.
\end{align}
Then, for $\ell = 2,\ldots ,L$, the equations~\eqref{eq:hidden-layer-first-node}--\eqref{eq:hidden-layer} become
\begin{align}\label{eq:dnn-matrix-form-1}
	x^\ell = A x^\ell + \begin{pmatrix}e^{-\alpha\theta}x^{\ell -1}_N\\0\\ \vdots \\ 0\end{pmatrix} + \alpha^{-1}(1-e^{-\alpha \theta})f(W^\ell x^{\ell -1}).
\end{align}
where $f$ is applied component-wise. By subtracting $A x^\ell$ from both sides of Eq.~\eqref{eq:dnn-matrix-form-1} and multiplication by the matrix
\begin{align}
	E := (\mathrm{Id} - A)^{-1} = \begin{pmatrix}
	1 & 0 & \cdots & \cdots & 0 \\
	e^{-\alpha\theta} & 1 & \ddots & & \vdots \\
	e^{-2\alpha\theta} & \ddots & \ddots & \ddots & \vdots \\
	\vdots & \ddots & \ddots & \ddots & 0 \\
	e^{-(N-1)\alpha\theta} & \cdots & e^{-2\alpha\theta} & e^{-\alpha\theta} & 1
	\end{pmatrix},
\end{align}
we obtain a matrix equation describing the $\ell$-th hidden layer
\begin{align}
	x^\ell 
	= \begin{pmatrix}e^{-\alpha\theta}x^{\ell -1}_N\\e^{-2\alpha\theta}x^{\ell -1}_N\\ \vdots \\ e^{-N\alpha\theta}x^{\ell -1}_N\end{pmatrix} +  \alpha^{-1}(1-e^{-\alpha \theta}) E f(W^\ell x^{\ell -1}).
\end{align}

The neural network~\eqref{eq:first-hidden-layer-first-node}--\eqref{eq:activation-output} obtained from delay system~\eqref{eq:semilinear-system}--\eqref{eq:activation-signal} can be trained by gradient descent~\cite{Stelzer2020}. The training parameters are the entries of the matrices $W^\mathrm{in}$ and $W^\mathrm{out}$, the step heights of the modulation functions $\mathcal{M}_d(t)$, and the bias signal $b(t)$.

\subsection{Network for large node distance $\theta$}
\label{subsec:largetheta}

In contrast to the general system~\eqref{eq:general-system}, the semilinear system~\eqref{eq:semilinear-system} with activation signal~\eqref{eq:activation-signal} does not only emulate a network of nodes for small distance $\theta$. It is also possible to choose large $\theta$. In this case, we can approximate the nodes given by Eq.~\eqref{eq:nodes} by the map limit
\begin{align}\label{eq:map-limit}
	 & x^\ell  = \alpha^{-1} f(a^\ell),  \\
	 & \text{where} \quad a^\ell  = W^\ell x^{\ell-1} \  \text{for} \  \ell>1 \quad  \text{and} \quad  a^1 =  g(W^\mathrm{in} u),
\end{align}
up to exponentially small terms. 

The reason for this limit behavior lies in the nature of the local couplings. Considering Eq.~\eqref{eq:semilinear-system}, one can interpret the parameter $\alpha$ as a time scale of the system, which determines how fast information about the system state at a certain time point decays while the system is evolving.
This phenomenon is related to the so-called instantaneous Lyapunov exponent \cite{Heiligenthal2011,Heiligenthal2013,Kinzel2013}, which equals $-\alpha$ in this case. As a result,
the local coupling between neighboring nodes emerges when only a small amount of time $\theta$ passes between the nodes. 
Hence, increasing $\theta$ one can reduce the local coupling strength until it vanishes up to a negligibly small value.
For a rigorous derivation of Eq.~\eqref{eq:map-limit}, we refer to \cite{Stelzer2020}.

The apparent advantage of the map limit case is that the obtained network matches a classical multilayer perceptron. Hence, known methods such as gradient descent training via the classical back-propagation algorithm \cite{Rumelhart1986} can be applied to the delay-induced network~\cite{Stelzer2020}.

The downside of choosing large values for the node separation $\theta$ is that the overall processing time of the system scales linearly with $\theta$. We need a period of time $T=N\theta$ to process one hidden layer. Hence, processing a whole network with $L$ hidden layers requires the time period $LT=LN\theta$. For this reason, the work~\cite{Stelzer2020} provides a modified back-propagation algorithm for small node separations to enable gradient descent training of networks with significant local coupling.

\section{Conclusions}

We have shown how networks of coupled maps with arbitrary topology and arbitrary size can be emulated by a single (possibly even scalar) DDE with multiple delays. Importantly, the coupling weights can be adjusted by changing  the modulations of the feedback signals. The network topology is determined by the choice of time-delays. As shown previously \cite{Appeltant2011,Larger2012,Brunner2013,Larger2017,Stelzer2020}, special cases of such networks are successfully applied for reservoir computing or deep learning. 

As an interesting conclusion, it follows that the temporal dynamics of DDEs can unfold arbitrary spatial complexity, which, in our case, is reflected by the topology of the unfolded network. In this respect, we shall mention previously reported spatio-temporal properties of DDEs \cite{Arecchi1992,Giacomelli1994,Giacomelli1996,Bestehorn2000,Giacomelli2012,Yanchuk2014,Yanchuk2015a,Yanchuk2015b,Kashchenko2016,Yanchuk2017}. These results show how in some limits, mainly for large delays, the DDEs can be approximated by partial differential equations.

Further, we remark that similar procedures have been used for deriving networks from systems of ordinary differential equations \cite{Haber2017,Lu2018,Chen2018}. However, in their approach, one should use an $N$-dimensional system of equations for implementing layers with $N$ nodes. This is in contrast to the DDE case, where the construction is possible with just a single-variable equation.

As a possible extension, a realization of adaptive networks using a single node with delayed feedback would be an interesting open problem. In fact, the application to deep neural networks in \cite{Stelzer2020} realizes an adaptive mechanism for the adjustment of the coupling weights. However, this adaptive mechanism is specially tailored for DNN problems. Another possibility would be to emulate networks with dynamical adaptivity of connections \cite{Berner}. 
The presented scheme can also be extended by employing delay differential-algebraic equations \cite{Mehrmann,Unger}.

\section*{Acknowledgements}

This work was funded by the ``Deutsche Forschungsgemeinschaft'' (DFG) in the framework of the project 411803875 (S.Y.) and IRTG 1740 (F.S.).

\end{document}